\documentclass[12pt,reqno]{article}
\usepackage{amsmath, amssymb}

\newtheorem{theorem}{Theorem}[section]
\newtheorem{definition}[theorem]{Definition}

\newtheorem{corollary}[theorem]{Corollary}
\newtheorem{lemma}[theorem]{Lemma}

\def \proof {\noindent {\bf Proof.}\ \ }

\def \endproof {{\mbox{}\nolinebreak\hfill\rule{2mm}{2mm}\par\medbreak}
}

\def \R {\mathbb{R}}

\def \E {\mathbb{E}}
\def \P {\mathbb{P}}

\def \one {{\bf 1}}
\def \PP {\mathcal{P}}
\def \SS {\mathcal{S}}
\def \a {\alpha}

\def \g {\gamma}
\def \e {\varepsilon}
\def \eps {\varepsilon}
\def \d {\delta}

\def \l {\lambda}

\def \s {\sigma}

\def \< {\langle}
\def \> {\rangle}

\def \diam {{\rm diam}}

\def \conv {{\rm conv}}

\def \vol {{\rm vol}}

\def \bi {B\bigl(L_\infty(\Omega)\bigr)}
\def \VC {{\rm VC}}

\begin{document}
\title {Entropy, dimension and the Elton-Pajor Theorem}
\author {S. Mendelson\footnote{Research School of
Information Sciences and Engineering, The Australian National
University, Canberra, ACT 0200, Australia, e-mail:
shahar.mendelson@anu.edu.au} \and
   R. Vershynin\footnote{
   Department of Mathematical Sciences,
   University of Alberta,
   Edmonton, Alberta T6G 2G1, Canada,
   e-mail: vershynin@yahoo.com}}
\date{December 17, 2001}

\maketitle
\begin{abstract}
The Vapnik-Chervonenkis dimension of a set $K$ in $\R^n$ is the
maximal dimension of the coordinate cube of a given size, which
can be found in coordinate projections of $K$. We show that the
VC dimension of a convex body governs its entropy. This has a
number of consequences, including the optimal Elton's theorem and
a uniform central limit theorem in the real valued case.
\end{abstract}


\section{Introduction}

Let $x_1, \ldots, x_n$ be vectors in the unit ball of a Banach
space, and assume that $\E \| \sum_{i=1}^n \e_i x_i \|  \ge  \d n$
for some number $\d > 0$, where $\e_1, \ldots, \e_n$ denote
independent Bernoulli random variables (taking values $1$ and
$-1$ with probability $1/2$). In 1983, J.~Elton \cite{E} proved an
important result that there exists a subset $\s$ of $\{1, \ldots,
n\}$ of size proportional to $n$ such that the set of vectors
$(x_i)_{i \in \s}$ is well equivalent to the $\ell_1$ unit-vector
basis. Specifically, there exist numbers $s, t > 0$, depending
only on $\d$, such that $|\s| \ge s n$ and $\| \sum_{i \in \s}
a_i x_i \|  \ge  t \sum_{i \in \s} |a_i|$ for all real numbers
$(a_i)$. This result was extended to the complex case by A.~Pajor
\cite{Pa}.

Several steps have been made towards finding asymptotically the
largest possible $s$ and $t$ in Elton's Theorem (\cite{Pa},
\cite{T}). Trivial upper bounds are that $s \le \d^2$, which
follows from the example of identical vectors, and $t \le \d$ as
demonstrated by shrinking the usual $\ell_1^n$ unit-vector basis.
One of the aims of this paper is to prove Elton's Theorem with $s
\ge c \d^2$ and $t \ge c \d$, where $c> 0$ is an absolute
constant. Furthermore, we show that $s$ and $t$ satisfy $\sqrt{s}
t \log^{2.1} (2 / t) \ge c \d$, which, as an easy example shows,
is optimal for all $\d$ up to a logarithmic factor. This improves
the result of M.~Talagrand from \cite{T}.

This theorem follows from new entropy estimates of a convex body
$K \subset [-1,1]^n  =  B_\infty^n$. We show that the entropy of
$K$ is controlled by its Vapnik-Chervonenkis dimension. This
parameter, denoted by $\VC(K,t)$, is defined for every $0 < t <
1$ as the maximal size of a subset $\s$ of $\{1, \ldots, n\}$,
such that the coordinate projection of $K$ onto $\R^\s$ contains
a coordinate cube of the form $x + [0, t]^\s$. This notion
carries over to convexity the ``classical'' concept of the VC
dimension, denoted by $\VC(A)$, and defined for subsets $A$ of
the discrete cube $\{0, 1\}^n$ as the maximal size of the subset
$\s$ of $\{1, \ldots, n\}$ such that $P_\s A  = \{0, 1\}^\s$,
where $P_\s$ is the coordinate projection onto the coordinates in
$\s$ (see \cite{L-T} \S 14.3).

Consider the unit ball $B_p^n$ of $\ell_p^n$, $1 \le p \le
\infty$, and let us look at the covering numbers $N(K,
n^{1/p}B_p^n, t)$, which are the minimal number of translates of
$t n^{1/p}B_p^n$ in $\R^n$ needed to cover $K$. A volumetric
bound on the entropy (which is the logarithm of the covering
numbers) shows that
$$
\log N(K, n^{1/p} B_p^n, t)  \le  \log(5 / t) \cdot n.
$$
One question is whether it is possible to replace the dimension
$n$ on the right-hand side of this estimate by the VC dimension
$\VC(K, ct)$, which is generally smaller? This is perfectly true
for the {\em Boolean} cube: the known theorem of R.~Dudley that
lead to a characterization of the uniform central limit property
in the Boolean case states that if $A  \subset \{0, 1\}^n$ then
$$
\log N(A, n^{1/2} B_2^n, t)  \le  C \log(2 / t) \cdot \VC(A).
$$
This estimate follows by a random choice of coordinates and an
application of the Sauer-Shelah Lemma (see \cite{L-T} Theorem
14.12). The same problem for convex bodies is considerably more
difficult, as to bound $\VC(K, t)$ one needs to find a cube in
$P_\s K$ with well separated faces, not merely disjoint. We prove
the following theorem.

\begin{theorem}                             \label{bp binfty}
   There are absolute constants $C, c > 0$ such that for every convex
   body $K \subset B_\infty^n$,
   every $1 < p < \infty$ and any $0 < t < 1$,
   \begin{equation} \label{bpn}
   \log N(K, n^{1/p} B_p^n, t)
     \le  Cp^2 \log^2 (2 / t) \cdot \VC(K, c t).
   \end{equation}
   Moreover,
   \begin{equation} \label{binfty}
   \log N(K, B_\infty^n, t)
     \le  C M^2 \log^2 (2 / t) \cdot \VC(K, c t),
   \end{equation}
   provided that either the right or the left hand side
   of \eqref{binfty} is larger than $t^M n$.
\end{theorem}

Let us comment on estimate \eqref{binfty}, which improves the
main lemma of \cite{A-B-C-H}. This bound can not hold in general
if the coefficient in front of the VC dimension depends only on
$t$ and not on $n$, since for $K = B_1^n$ we have $\VC(K, t) = 2
/ t$ and $\log N(K, B_\infty^n, t) \ge \log n$. Next,
\eqref{binfty} is best complemented by the easy lower bound
$$
\log N(K, B_\infty^n, t)  \ge  \VC(K, c t),
$$
for some absolute constant $c > 0$, which follows from the
definition of the VC dimension and by a comparison of volumes.
These two bounds show that the $\|\cdot\|_\infty$-entropy of $K$
is governed by the VC dimension of $K$, up to a logarithmic
factor in $t$.

The relation to the Elton-Pajor Theorem is the following. If $K$
is a symmetric convex body, then $\VC(K, t)$ is the maximal
cardinality of a subset $\s$ of $\{1, \ldots, n\}$ such that $\|
\sum_{i \in \s} a_i e_i \|_{K^\circ} \ge  (t / 2) \sum_{i \in \s}
|a_i|$ for all real numbers $(a_i)$, where $e_i$ are the
canonical unit vectors in $\R^n$ and $K^\circ$ is the polar of
$K$. Note that if $(g_i)$ are independent standard gaussian random
variables then $E\|\sum_{i=1}^n \e_i e_i\| \leq 2E\|\sum_{i=1}^n
g_i e_i\|$ for every norm (\cite{L-T} \S 4.5). Therefore, our
problem reduces to finding a bound on
$$
E = \E \Big\| \sum_{i=1}^n g_i e_i \Big\|_{K^\circ}
$$
in terms of the VC-dimension of $K$. The latter is relatively easy
once we know \eqref{bpn}. Indeed, replacing the entropy by the VC
dimension in Dudley's entropy inequality it follows that there
are absolute constants $C$ and $c$ such that
\begin{equation}                                \label{E}
E \le  C  \int_{c E / \sqrt{n}}^{\infty}
              \sqrt{ \log N(K, B_2^n, t) } \; d t
    \le  C \sqrt{n} \int_{c E / n}^1  \sqrt{\VC(K, ct)} \log(2 / t) \;
dt.
\end{equation}
This inequality improves of the main theorem of M.~Talagrand in
\cite{T}. Elton's Theorem with optimal asymptotics follows from
\eqref{E} by comparing the integrand to an appropriately chosen
integrable function.

We present a few other applications to convexity. Inequality
\eqref{E} can be applied, as in \cite{T}, to compare two
geometric properties of a Banach space called type and infratype.
Recall that a Banach space $X$ is of gaussian type $p$ if there
exists some $M > 0$ such that for all $n$ and all sequences of
vectors $(x_i)_{i \le n}$,
\begin{equation}                            \label{Tp}
   \E \Big\| \sum_{i = 1}^n g_i x_i \Big\| \le  M \Big( \sum_{i =
   1}^n \|x_i\|^p \Big)^{1/p}.
\end{equation}
The best possible constant $M$ in this inequality is denoted by
$T_p(X)$. Next, $X$ has infratype $p$ if there exists some $M >
0$ such that for all $n$ and all sequences of vectors $(x_i)_{i
\le n}$, we have
\begin{equation}                            \label{Ip}
   \min_{\eta_i = \pm 1} \Big\| \sum_{i = 1}^n \eta_i x_i \Big\|
   \le  M \Big( \sum_{i = 1}^n \|x_i\|^p \Big)^{1/p}.
\end{equation}
The best possible constant $M$ in this inequality is denoted by
$I_p(X)$.

M.~Talagrand proved in \cite{T} that if $1 < p < 2$ then $T_p(X)
\le C(p) I_p(X)^2$, where $C(p)$ is a constant which depends only
on $p$. It is not known whether the square can be removed.
Moreover, the situation for $p = 2$ is unknown in general, but
\eqref{E} can be used to show that there is an absolute constant
$C$ such that for any $n$ dimensional Banach space $X$,
$$
T_2 (X)  \le  I_2 (X)  \cdot  C \log^2 \Big( \frac{n}{I_2(X)^2}
\Big)
   \le  I_2 (X)  \cdot  C \log^2 n.
$$

Finally, we present an application of Theorem \ref{bp binfty} to
empirical processes. We use a version of \eqref{bpn} to bound the
entropy of an arbitrary subset of $B_\infty^n$ using a
scale-sensitive version of the ``classical" VC dimension, known
as the fat-shattering dimension. In particular we show that if
$F$ is a class of uniformly bounded functions, which has a
relatively small fat-shattering dimension, then it satisfies the
uniform central limit theorem for any probability measure. This
extends Dudley's characterization for VC classes to the
real-valued case.

The paper is organized as follows. In Section \ref{s: bpn} we
prove the bound for the $B_p^n$-entropy in abstract finite
product spaces, and then derive \eqref{bpn} by approximation.
Actually, the convexity of $K$ plays a very little role in these
results, and similar entropy bounds hold for arbitrary susets of
$B_\infty^n$. In Section \ref{s: binfty} we prove \eqref{binfty}
for the $B_\infty^n$-entropy by reducing it to \eqref{bpn}
through an independent lemma that compares the $B_p^n$-entropy to
the $B_\infty^n$-entropy. In Section \ref{s: convex} we apply
\eqref{bpn} to convex bodies. In particular, we deduce Elton's
Theorem and the infratype results.  Finally, in Section \ref{s:
empirical measures} we apply \eqref{bpn} to empirical processes.

Throughout this article, positive absolute constants are denoted
by $C$ and $c$. Their values may change from line to line, or
even within the same line.

ACKNOWLEDGEMENTS: The second author is thankful to Mark Rudelson
who contributed a lot of effort and enthusiasm to the paper.
Warmest thanks are to Nicole Tomczak-Jaegermann for her constant
support. The second author also acknozledges a support from the
Pacific Institute of Mathematical Sciences, and thanks the
Department of Mathematicql Sciences of the University of Alberta
for hospitality.

\section{$B_p^n$-entropy in abstract product spaces}        \label{s:
bpn}

We will introduce and work with the notion of the VC dimension in
an abstract setting that encompasses both classes considered in
the introduction, the subsets of the discrete cube $\{0, 1\}^n$
and the class of convex bodies in $\R^n$.

We call a map $d : T \times T \to \R_+$ a {\em quasi-metric} if
$d$ is symmetric and reflexive (that is, $\forall x,y$, $d(x,y) =
d(y,x)$ and $d(x,x) = 0$). We say that points $x$ and $y$ in $T$
are separated if $d(x,y)>0$. Thus, $d$ does not necessarily
separate points or satisfy the triangle inequality.

\begin{definition}
   Let $(T,d)$ be a quasi-metric space
   and let $n$ be a positive integer.
   For a set $A  \subset T^n$ and $t  >  0$,
   the VC-dimension $\VC(A, t)$ is the maximal cardinality
   of a subset $\s  \subset \{1, \ldots, n \}$ such that
   the inclusion
   \begin{equation}
     P_\s A  \supseteq  \prod_{i \in \s} \{a_i, b_i\}
\label{inclusion}
   \end{equation}
   holds for some points $a_i, b_i  \in T$, $ i \in \s$
   with $d(a_i, b_i)  \ge  \d$.
   If no such $\s$ exists, we set $\VC(A, t) = 0$.
   When there is a need to specify the underlying metric,
   we denote the VC dimension by $\VC_d (A,  t)$.
\end{definition}

Since $\VC(A, t)$ is decreasing in $t$ and is bounded by $n$,
which is the ``usual" dimension of the product space, the limit
$$
\VC(A)  :=  \lim_{t \to 0_+} \VC(A, t)
$$
always exists. Equivalently, $\VC(A)$ is the maximal cardinality
of a subset $\s  \subset \{1, \ldots, n \}$ such that
\eqref{inclusion} holds for some pairs $(a_i, b_i)$ of separated
points in $T$.

This definition is an extension of the ``classical" VC dimension
for subsets of the discrete cube $\{0, 1\}^n$, where we think of
$\{0,1\}$ as a metric space with the $0 -1$ metric. Clearly, for
any set $A \subset \{0, 1\}^n$ the quantity $\VC(A, t)$ does not
depend on $0 < t < 1$, and hence
$$
\VC(A)  =  \max \Big\{ |\s| : \;  \s \subset \{1, \ldots, n\}, \;
   P_\s A = \{0, 1\}^\s  \Big\},
$$
which is precisely the ``classical" definition of the VC
dimension.

The other example discussed in the introduction was the VC
dimension of convex bodies. Here $T = \R$ or, more frequently, $T
= [-1,1]$, both with respect to the usual metric. If $K \subset
T^n$ is a convex body, then $\VC(K, t)$ is the maximal
cardinality of a subset $\s \subset \{1, \ldots, n \}$ for which
the inclusion
$$
P_\s K  \supseteq  x + (t / 2) B_\infty^\s
$$
holds for some vector $x \in \R^\s$ (which automatically lies in
$P_\s K$). It is easy to see that if $K$ is symmetric, we can set
$x = 0$. Also note that for every convex body $\VC(K) = n$.

The main results of this article rely on (and are easily reduced
to) a discrete problem: to estimate the VC-dimension of a set in a
product space $T^n$, where $(T, d)$ is a {\em finite}
quasi-metric space. $T^n$ is usually endowed with the normalized
Hamming quasi-metric $d_n (x, y)  =  n^{-1} \sum_{i = 1}^n
d(x(i), y(i))$ for $x, y \in T^n$.

In the main result of this section we bound the entropy of a set
$A \subset T^n$ with respect to $d_n$ in terms of $\VC(A)$.

\begin{theorem}                           \label{in product}
   Let $(T, d)$ be a finite quasi-metric space
   with $\diam (T)  \le 1$,
   and set $n$ to be a positive integer.
   Then, for every set $A  \subset T^n$
   and every $0 <  \e  <  1$,
   $$
   \log N(A, d_n, \e)  \le  C \log^2 (|T| / \e) \cdot \VC(A),
   $$
   where $C$ is an absolute constant.
\end{theorem}

Before presenting the proof, let us make two standard
observations. We say that points $x, y \in T^n$ are separated on
the coordinate $i_0$ if $x(i_0)$ and $y(i_0)$ are separated.
Points $x$ and $y$ are called $\e$-separated if $d_n(x,y) \geq
\e$.

Clearly, if $A'$ is a maximal $\e$-separated subset of $A$ then
$|A'| \geq N(A, d_n, \e)$. Moreover, the definition of $d_n$ and
the fact that $\diam(T) \leq 1$ imply that every two distinct
points in $A'$ are separated on at least $\e n$ coordinates. This
shows that Theorem \ref{in product} may be reduced to the
following statement.

\begin{theorem}                             \label{thm:newproduct}
   Let $(T,d)$ be a quasi-metric space
   for which $\diam (T) \leq 1$.
   Let $0 < \eps < 1$ and
   consider a set $A \subset T^n$ such that
   every two distinct points in $A$ are separated
   on at least $\e n$ coordinates.
   Then
   \begin{equation}                      \label{VC large wanted}
     \log |A|  \le  C \log^2 (|T| / \e) \cdot \VC(A).
   \end{equation}
\end{theorem}

The first step in the proof of Theorem \ref{thm:newproduct} is a
probabilistic extraction principle, which allows one to reduce
the number of coordinates without changing the separation
assumption by much. Its proof is based on a simple discrepancy
bound for a set system.

\begin{lemma}                           \label{discrepancy}
   There exists an absolute constant $c > 0$ for which
   the following holds.
   Let $\e > 0$ and assume that $\SS$ is a system of subsets
   of $\{1, \ldots, n\}$ which satisfies that each $S \in \SS$ contains
   at least $\e n$ elements.
   Let $k \le n$ be an integer such that
   $\log |\SS|   \le  c \e k$.
   Then there exists a subset $I \subset \{1, \ldots, n\}$ of
   cardinality $|I| = k$,
   such that
   $$
   |I \cap S|  \ge   \e k/4
   \ \ \ \text{for all $S \in \SS$}.
   $$
\end{lemma}

\proof If $|\SS| = 1$ the lemma is trivially true, hence we may
assume that $|\SS| \ge 2$. Let $0 < \d < 1/2$ and set $\d_1,
\ldots, \d_n$ to be $\{0,1\}$-valued independent random variables
with $\E \d_i = \d$ for all $i$. By the classical bounds on the
tails of the binomial law (see \cite{H}, or \cite{L-T} 6.3 for
more general inequalities), there is an absolute constant $c_0 >
0$ for which
\begin{equation}                        \label{hoeffding}
   \P \Big\{  \big| \sum_{i =1}^n (\d_i - \d)  \big| >  \frac{1}{2}
   \d n \Big\}
   \le  2 \exp (- c_0 \d n).
\end{equation}
Let $\d = {k}/{2n}$ and consider the random set $I = \{ i: \; \d_i
= 1 \}$. For any set $B  \subset  \{1, \ldots, n \}$, $|I \cap
B|  = \sum_{i \in B} \d_i$. Then \eqref{hoeffding} implies that
$$
\P \{  |I \cap B| \ge  \d |B|/2 \}
   \ge  1 - 2 \exp (- c_0 \d |B|).
$$
Since for every $S \in \SS$, $|S| > \e n$, then
$$
\P \{  |I \cap S|  \ge  \e k/4 \}
   \ge  1 - 2 \exp (- \frac{1}{2} c_0 \e k).
$$
Therefore,
$$
\P \Big\{  \forall S \in \SS, \; |I \cap S|  \ge  \frac{1}{4} \e k
\Big\}
   \ge  1 - 2 |\SS| \exp (-  \frac{1}{2}  c_0 \e k).
$$
By the assumption on $k$, this quantity is larger than $1/2$ (with
an appropriately chosen absolute constant $c$). Moreover, by a
similar argument, $|I| \le k$ with probability larger than $1/2$.
This proves the existence of a set $I$ satisfying the assumptions
of the lemma.
\endproof

\noindent{\bf Proof of Theorem \ref{thm:newproduct}. } We may
assume that $|T| \ge 2$, $\e  \le  1/2$, $n \ge 2$ and $\max(4,
\exp(4 c))  \le  |A|  \le |T|^n$, where $0 < c < 1$ is the
constant in Lemma \ref{discrepancy}. The first step in the proof
is to use previous lemma, which enables one to make the
additional assumption that $\log |A|  \ge  c \e n/4$. Indeed,
assume that the converse inequality holds, and for every pair of
distinct points $x, y \in A$, let $S(x, y) \subset \{1, \ldots,
n\}$ be the set of coordinates on which $x$ and $y$ are
separated. Put $\SS$ to be the collection of the sets $S(x,y)$
and let $k$ be the minimal positive integer for which $\log|\SS|
\le  c \e k$. Since $|A| \leq |\SS| \le |A|^2$, then
$$
c \e (k-1)  \le  \log |\SS| \le 2 \log |A|  \le \frac{1}{2} c \e
n,
$$
which implies that $1 \le k \le n$. Thus, by Lemma
\ref{discrepancy} there is a set $I \subset \{1, \ldots, n\}$,
$|I| = k$, with the property that every pair of distinct points
$x, y \in A$ is separated on at least $ \e |I|/4$ coordinates in
$I$. Also, since $4c \leq \log |A| \leq \log |\SS| \leq c\e k$,
then $\eps |I|/4 \geq 1$ and thus $|P_I A|=|A|$. Clearly, to
prove the assertion of the theorem for the set $A \subset T^n$,
it is sufficient to prove it for the set $P_I A \subset T^I$ (with
$|I|$ instead of $n$), whose cardinality already satisfies $\log
|P_I A| =  \log |A| \ge  c \e (k - 1)/2
   \ge   c \e |I|/4$.
Therefore, we can assume that $|A| = \exp(\a n)$ with $\a > c \e$
for some absolute constant $c$.

The next step in the proof is a counting argument, which is based
on the proof of Lemma 3.3 in \cite{A-B-C-H} (see also \cite{B-L}).

A set is called a {\em cube} if it is of the form $D_\s =
\prod_{i \in \s} \{ a_i, b_i \}$, where $\s$ is a subset of $\{1,
\ldots, n\}$ and $a_i, b_i  \in  T$. We will be interested only in
{\em large cubes}, which are the cubes in which $a_i$ and $b_i$
are separated for all $i \in \s$. Given a set $B \subset T^n$, we
say that a cube $D_\s$ {\em embeds} into $B$ if $D_\s  \subset
P_\s B$. Note that if a large cube $D_\s$ with $|\s| \ge v$ embeds
into $B$ then $\VC(B) \ge v$.

For all $m \ge 2$, $n \ge 1$ and $0 < \e \le 1/2$, let $t_\e (m,
n)$ denote the maximal number $t$ such that for every set $B
\subset T^n$, $|B| = m$, which satisfies the separation condition
we imposed (that is, every distinct points $x,y \in B$ are
separated on at least $\eps n$ coordinates), there exist $t$
large cubes that embed into $B$. If no such $B$ exists, we set
$t_\e (m, n)$ to be infinite.

The number of possible large cubes $D_\s$ for $|\s| \le v$ is
smaller than $\sum_{k=1}^v  \binom{n}{k} |T|^{2 k}$, as for every
$\s$ of cardinality $k$ there are less than $|T|^{2 k}$
possibilities to choose $D_\s$. Therefore, if $t_\e (|A|, n) \ge
\sum_{k=1}^v \binom{n}{k} |T|^{2 k}$, there exists a large cube
$D_\s$ for some $|\s| \geq v$ that embeds into $A$, implying that
$\VC(A) \ge  v$. Thus, to prove the theorem, it suffices to
estimate $t_\e (m, n)$ from below. To that end, we will show that
for every $n \geq 2$, $m \geq 1$ and $0<\e \leq 1/2$,
\begin{equation} \label{eq:claim}
t_\e ( 2 m \cdot |T|^2 / \e, n)  \ge  2 t_\e (2 m, n - 1).
\end{equation}
Indeed, fix any set $B \subset T^n$ of cardinality $|B| = 2 m
\cdot |T|^2 / \e$, which satisfies the separation condition
above. If no such $B$ exists then $t_\e ( 2 m \cdot |T|^2 / \e,
n) = \infty$, and (\ref{eq:claim}) holds trivially. Split $B$
arbitrarily into $m \cdot |T|^2 / \e$ pairs, and denote the set
of the pairs by $\PP$. For each pair $(x, y) \in \PP$ let $I
(x,y) \subset \{1, \ldots, n\}$ be the set of the coordinates on
which $x$ and $y$ are separated, and note that by the separation
condition, $|I (x, y)| \ge \e n$.

Let $i_0$ be the random coordinate, that is, a random variable
uniformly distributed in $\{1, \ldots, n\}$. The expected number
of the pairs $(x, y) \in \PP$ for which $i_0 \in  I(x,y)$ is
$$
\E \sum_{(x, y) \in \PP} \one_{ \{i_0  \in  I(x,y) \} } =
\sum_{(x, y) \in \PP}  \P \{i_0  \in  I(x,y) \} \ge  |\PP| \cdot
\e =  m |T|^2.
$$
Hence, there is a coordinate $i_0$ on which at least $m |T|^2$
pairs $(x, y) \in \PP$ are separated. By the pigeonhole principle,
there are at least $m |T|^2 / \binom{|T|}{2} \ge  2 m$ pairs
$(x,y) \in \PP$ for which the (unordered) set $\{ x(i_0), y(i_0)
\}$ is the same.

Let $I = \{ 1, \ldots, n \} \setminus \{ i_0\}$. It follows that
there are two subsets of $B$, denoted by $B_1$ and $B_2$, such
that $|B_1| = |B_2| = 2 m$ and
$$
B_1  \subset  \{b_1\}  \times  T^I,  \; \; B_2  \subset  \{b_2\}
\times  T^I
$$
for some separated points $b_1, b_2 \in T$. Clearly, the set
$B_1$ satisfies the separation condition and so does $B_2$. It is
also clear that if a large cube $D_\s$ embeds into $B_1$, then it
also embeds into $B$, and the same holds for $B_2$. Moreover, if
the same cube $D_\s$ embeds into both $B_1$ and $B_2$, then the
large cube $\{b_1, b_2\} \times D_\s$ embeds into $B$ (since
$\{b_1, b_2\} \times D_\s  \subset P_{\{i_0\} \cup \s} B$).
Therefore, $t_\e (|B|, n)  \ge  2 t_\frac{\e n}{n-1} (|B_1|, n -
1)  \ge  2 t_\e (|B_1|, n - 1)$, establishing \eqref{eq:claim}.

\qquad

Since $t_\e (2, n)  \ge  1$, an induction argument yields that
$t_\e ( 2 (|T|^2 / \e)^r, n )  \ge  2^r$ for every $r \ge 1$.
Thus, for every $m \ge 4$
$$
t_\e (m, n)  \ge  m^{\frac{1}{2 \log(|T|^2 / \e)}}.
$$
(It is remarkable that the right hand side does not depend on
$n$). Therefore, $\VC(A) \ge v$ provided that $v$ satisfies
\begin{equation}                \label{requirement on d}
   t_\e (|A|, n)  \ge  \exp \Big( \frac{\a n}{2 \log(|T|^2 / \e)} \Big)
   \ge  \sum_{k=1}^v  \binom{n}{k} |T|^{2 k}.
\end{equation}
To estimate $v$, one can bound the right-hand side of
\eqref{requirement on d} using Stirling's approximation $\sum_{k
= 1}^v \binom{n}{k}  \le   [ \g^\g (1 - \g)^{1 - \g} ]^{-n}$,
where $\g = v / n  \le  1/2$. It follows that for $v  \le  n/2$,
$\sum_{k = 1}^v \binom{n}{k} |T|^{2 k} \le  ( \frac{|T| n}{v}
)^{2 v}$. Taking logarithms in \eqref{requirement on d}, we seek
integers $v \le n/2$ satisfying that
$$
\frac{\a n}{2 \log(|T|^2 / \e)}  \ge  2 v \log \Big( \frac{|T|
n}{v} \Big).
$$
This holds if
$$
v  \le  \cdot \left( \frac{\a n}{\log(|T|^2 / \e)} \right)
             \Big/ 8 \log \left( \frac{4 |T| \log(|T|^2 / \e)}{ \a}
\right),
$$
proving our assertion since $\a > c \e$.
\endproof

\begin{corollary}                                \label{in lattice}
   Let $n \ge 2$ and $p \ge 2$ be integers, set
   $0 < \e  < 1$ and $q > 0$.
   Consider a set $A \subset \{ 1, \ldots, p \}^n$
   such that for every two distinct points $x,y \in A$,
   $|x(i) - y(i)| \ge q$ for at least $\e n$ coordinates $i$.
   Then
   $$
   \log |A|  \le  C \log^2 (p / \e) \cdot \VC(A, q).
   $$
\end{corollary}

\proof We can assume that $q  \ge  1$. Define the following
quasi-metric on $T = \{1, \ldots, p\}$:
$$
d(a, b)  =
\begin{cases}
   0  &  \text{if $|a - b| < q$}, \\
   1  &  \text{otherwise.}
\end{cases}
$$
Then $N (A, d_n, \e)  =  |A|$. By Theorem \ref{in product},
$$
\log |A|  \le  C \log^2 (p / \e) \cdot \VC_d (A),
$$
which completes the proof by the definition of the metric $d$.
\endproof

Now we pass from the discrete setting to the ``continuous" one -
namely, we study subsets of $B_\infty^n$. Recall that the
Minkowski sum of two convex bodies $A,B \subset \R^n$ is defined
as $A+B = \{a+b | \; a \in A, \ b \in B \}$.

\begin{corollary}                           \label{in binfty}
   For every $A \subset B_\infty^n$,
   $0 < t < 1$ and $0 < \e < 1$,
   $$
   \log N (A, \sqrt{n} B_2^n, t)
   \le  C \log^2 (2 / t \e)  \cdot  \VC(A + \e B_\infty^n, t / 2).
   $$
\end{corollary}

\proof Clearly, we may assume that $\e \le t/ 4$. Put $p =
\frac{1}{2 \e}$ and let
$$
T  =  \{ -2 \e p, -2 \e (p-1), \ldots, -2 \e, 0, 2 \e, \ldots, 2
\e (p-1), 2 \e p \}.
$$
Since $t- \e  >  3 t / 4$, then by approximation one can find a
subset $A_1  \subset  T^n$ for which $A_1  \subset  A + \e
B_\infty^n$ and $N(A_1, \sqrt{n} B_2^n, t - \e)  \ge  N (A,
\sqrt{n} B_2^n, t)$. Therefore, there exists a subset $A_2
\subset A_1$ of cardinality $|A_2| \ge  N (A, \sqrt{n} B_2^n,
t)$, which is $\frac{3 t}{4}\sqrt{n}$-separated with respect to
the $\|\cdot\|_2$-norm. Note that every two distinct points $x, y
\in A_2$ satisfy that
$$
\sum_{i = 1}^n |x(i) - y(i)|^2  \ge ({9 t^2}/{16}) n \geq {t^2}n/2
$$
and that $|x(i) - y(i)|^2 \le 4$ for all $i$. Hence $|x(i) - y(i)|
\ge {t}/{2}$ on at least ${t^2 n}/{16}$ coordinates $i$. By
Corollary \ref{in lattice} applied to $A_2$,
$$
\log |A_2|  \le  C \log^2 (2 / t \e) \cdot \VC(A_2, t / 2),
$$
and since $A_2  \subset  A_1 \subset A + \e B_\infty^n$, our
claim follows.
\endproof

 From this we derive the entropy estimate \eqref{bpn}.

\begin{corollary}                   \label{convex body}
   There exists an absolute constant $C$
   such that for any convex body
   $K \subset B_\infty^n$ and every $0<t<1$,
   $$
   \log N(K, \sqrt{n} B_2^n, t)
   \le  C \log^2 (2 / t) \cdot \VC(K, t / 4).
   $$
\end{corollary}

\proof This estimate follows from Corollary \ref{in binfty} by
selecting $\e = t / 4$ and recalling the fact that for every
convex body $K \subset \R^n$ and every $0 < b < a$,
$$
\VC(K + b B_\infty^n, a)  \le  \VC(K, a - b).
$$
The latter inequality is a consequence of the definition of the
VC-dimension and the observation that if $0 < b < a$ are such
that $a B_\infty^n  \subset K + b B_\infty^n$, then $(a - b)
B_\infty^n \subset K$.
\endproof

Note that Corollary \ref{in binfty} and Corollary \ref{convex
body} can be extended to the case where the covering numbers are
computed with respect to $n^{1/p}B_p^n$ for $1 < p < \infty$,
thus establishing the complete claim in \eqref{bpn}.

\section{$B_\infty^n$-entropy}          \label{s: binfty}

In this section we prove estimate \eqref{binfty}, which improves
the main combinatorial result in \cite{A-B-C-H}. Our result can
be equivalently stated as follows.

\begin{theorem}                                 \label{binfty thm}
   Let $K \subset B_\infty^n$ be a convex body, set $t > 0$ and put
   $v = \VC(K, t / 8)$. Then,
   \begin{equation}                          \label{NK small}
     \log N(K, B_\infty^n, t)  \le  C v \cdot \log^2 (n / t v),
   \end{equation}
   where $C$ is an absolute constant.
\end{theorem}

This estimate should be compared with the Sauer-Shelah lemma for
subsets of the Boolean cube $\{0,1\}^n$. It says that if $A
\subset \{0,1\}^n$ then for $v = \VC(K)$ we have $|A| \le
\binom{n}{0} + \binom{n}{1} + \ldots + \binom{n}{v}$, so that
$$
\log |A|  \le  2 v \cdot \log(n/v)
$$
(and note that, of course,  $|A| = N(K, B_\infty^n, t)$ for all
$0 < t < 1/2$).

We reduce the proof of \eqref{binfty thm} to an application of the
$B_p^n$-entropy estimate \eqref{bpn}. As a start, note that for $p
= \log n$, $B_\infty^n  \subset  n^{1/p} B_p^n  \subset  e
B_\infty^n$. Therefore, an application of \eqref{bpn} for this
value of $p$ yields
$$
\log N(K, B_\infty^n, t)  \le  C v \cdot \log^2 (n / t),
$$
which is slightly worse than \eqref{NK small}.

To deduce \eqref{NK small} we need a result that compares the
$B_\infty^n$-entropy to the $B_p^n$-entropy, and which may be
useful in other applications as well.

\begin{lemma}                           \label{separation}
   There is an absolute constant $c > 0$
   such that the following holds.
   Let $A$ be a subset of $B_\infty^n$ such that
   every two distinct points $x, y \in A$ satisfy
   $\|x - y\|_\infty  \ge  t$.
   Then, for every integer $1 \le k \le n/2$,
   there exists a subset $A'  \subset A$
   of cardinality
   $$
   |A'|  \ge  \binom{n}{k}^{-1}  (c t)^k |A|,
   $$
   with the property that every two distinct points in $A'$ satisfy 
that
   $|x(i) - y(i)|  \ge t / 2$
   for at least $k$ coordinates $i$.
\end{lemma}

\proof We can assume that $0 < t < 1/8$. Set $s = t / 2$. The
separation assumption imply that $N(A, B_\infty^n, s) \ge  |A|$.
Denote by $D_k$ the set of all points $x$ in $\R^n$ for which
$|x(i)| \ge 1$ on at most $k$ coordinates $i$. One can see that
$N(A,D_k, s) = N(A, s D_k, 1) = N(A, s D_k \cap 3 B_\infty^n,
1)$. Then, by the submultiplicative property of the covering
numbers,
\begin{eqnarray}                                \label{submult}
N(A, B_\infty^n, s)
   &\le&  N(A, s D_k \cap 3 B_\infty^n, 1)
              \cdot N \big( s D_k \cap 3 B_\infty^n, B_\infty^n, s 
\big)
\nonumber\\
    &\le&  N(A, s D_k, 1)
              \cdot N \big( s D_k \cap 3 B_\infty^n, B_\infty^n, s
\big).
\end{eqnarray}
To bound the second term, write $D_k$ as
$$
D_k  =  \bigcup_{|\s| = k}  \Big( \R^\s + (-1,1)^{\s^c} \Big),
$$
where the union is taken with respect to all subsets $\s \subset
\{1, \ldots, n\}$, and the sum in the right-hand side is the
Minkowski sum. Thus,
$$
s D_k  \cap 3 B_\infty^n
   =  \bigcup_{|\s| = k}  \Big( 3 B_\infty^\s + (-s, s)^{\s^c} \Big).
$$
Denote by $N' (A, B, t)$ the number of translates of $t B$ by
vectors in $A$ needed to cover $A$. Therefore,
\begin{eqnarray*}
N \big( s D_k \cap 3 B_\infty^n, B_\infty^n, s \big)
   &\le&  \sum_{|\s| = k}
              N \big( 3 B_\infty^\s + (-s, s)^{\s^c}, B_\infty^n, s
\big) \\
   &\le&  \sum_{|\s| = k}
              N' (3 B_\infty^\s, B_\infty^n, s).
\end{eqnarray*}
The latter inequality holds because any cover of $3 B_\infty^\s$
by translates of $s B_\infty^n$ automatically covers $3
B_\infty^\s + (-s, s)^{\s^c}$. Hence, for some absolute constant
$C$,
\begin{eqnarray*}
N \big( s D_k \cap 3 B_\infty^n, B_\infty^n, s \big)
   &\le&  \binom{n}{k} N' (3 B_\infty^k, B_\infty^k, s) \\
   &\le&  \binom{n}{k} (C / s)^k
\end{eqnarray*}
by a comparison of the volumes, and by \eqref{submult} we obtain
$$
N(A, D_k, s)  \ge  \binom{n}{k}^{-1} (c s)^k  N(A, B_\infty^n, s)
\ge  \binom{n}{k}^{-1} (c t)^k |A|,
$$
from which the statement of the lemma follows by the definition
of $D_k$.
\endproof

Now we can compare the $B_\infty^n$-entropy of $K$ to the $B_1^n$
entropy of $K$.

\begin{corollary}
   Let $A \subset B_\infty^n$ be a set,
   and set $0 < t < 1$ and $0 < \e < t/8$.
   Then
   $$
   N(A, B_\infty^n, t)
   \le  \Big( \frac{C}{\e} \Big)^{(2 \e / t) n}
         N(A, n B_1^n, \e),
   $$
   where $C$ is an absolute constant.
\end{corollary}

\proof Note that the set $A'$ in the conclusion of Lemma
\ref{separation} is such that every two distinct points $x, y \in
A'$ satisfy $\|x - y\|_1  \ge  (t / 2) k$. Thus $A'$ is $(t / 2)
k$-separated in the $\|\cdot\|_1$-norm, implying that $|A'|  \le
N(A, B_1^n, (t / 4) k)$. By Lemma \ref{separation},
$$
N(A, B_\infty^n, t)  \le  \binom{n}{k} (C / t)^k
    N(A, B_1^n, (t/4) k)
\le  \Big( \frac{C n}{t k} \Big)^{2 k}
    N \big( A, nB_1^n, \frac{t k}{4 n} \big).
$$
The conclusion follows by choosing $k$ which satisfies $\frac{t
k}{4 n} = \e$.
\endproof


\noindent {\bf Proof of Theorem \ref{binfty thm}. } Fix $0 < t <
1$, and let $\a$ be defined by $\log N(K, B_\infty^n, t)  =
\exp(\a n)$. Hence, there exists a set $A \subset K$ of
cardinality $|A| = \exp(\a n)$, where every two distinct points
$x, y \in A$ satisfy that $\|x - y\|_\infty  \ge  t$. Applying
Lemma \ref{separation} we obtain a subset $A' \subset A \subset
K$ of cardinality
$$
|A'|  \ge  \binom{n}{k}^{-1}  (c t)^k e^{\a n},
$$
such that for every two distinct points in $A'$, $|x(i) - y(i)|
\ge t / 2$ on at least $k$ coordinates $i$. Selecting $k = \frac{c
\a n}{\log(2 / t \a)}$ we see that $|A'|  \ge  e^{\a n/ 2}$.

The proof is completed by discretizing $A'$ and applying Corollary
\ref{in lattice} with $p = 4 / t$ and $\e = k / n$ in the same
manner as we did in the previous section. Therefore
\begin{eqnarray*}
\a n / 2 = \log |A'|
   &\le&  C \log^2 \big( \frac{4 n}{t k} \big)
             \cdot \VC (A' + (t / 4) B_\infty^n, t / 2) \\
   &\le&  C \log^2 (1 / t \a)  \cdot \VC(K, t / 4),
\end{eqnarray*}
and thus
$$
\a n  \le  c \log^2 (n / t v) \cdot v,
$$
as claimed.
\endproof

\section{Applications to convex bodies}             \label{s: convex}

We start by presenting an improvement of the main result of
M.~Talagrand from \cite{T}.
\begin{theorem} \label{thm:talagrand}
There are absolute constants $C, c > 0$ such that for every
convex body $K \subset B_\infty^n$
\begin{equation*}
E \leq C\sqrt{n} \int_{cE/n}^{1} \sqrt{\VC(K,ct)}\log (2/t) dt,
\end{equation*}
where $E=\E\|\sum_{i=1}^n g_i e_i\|_{K^\circ}$, and
$(e_i)_{i=1}^n$ is the canonical vector basis in $\R^n$.
\end{theorem}

For the proof, we need a few standard definitions and facts from
the local theory of Banach spaces, which may be found in
\cite{MS}.

Given an integer $n$, let $S^{n-1}$ be the unit Euclidean sphere
with the normalized Lebesgue measure $\s_n$, and for every
measurable set $A \subset \R^n$ denote by $\vol{A}$ its Lebesgue
measure in $\R^n$. For a convex body $K$ in $\R^n$, put $M_K =
\int_{S^{n-1}} \|x\|_K \; d \s_n(x)$ and let $M^*_K$ denote
$M_{K^\circ}$, where $K^\circ$ is the polar of $K$. Recall that
for any two convex bodies $K$ and $L$, $M^*_{K + L} \le M^*_K +
M^*_L$. Urysohn's inequality states that $\big(
\frac{\vol(K)}{\vol(B_2^n)} \big)^{1/n}  \le  M^*_K$.

Next, put $\ell(K)  =  \E \| \sum_{i=1}^n g_i e_i \|_K$, where
$(g_i)_{i=1}^n$ are independent standard gaussian random
variables and $(e_i)_{i=1}^n$ is the canonical basis of $\R^n$.
It is well known that $\ell(K)  =  c_n \sqrt{n} M_K$, where $c_n
< 1$ and $c_n \to 1$ as $n \to \infty$. Recall that by Dudley's
inequality (see \cite{Pi}) there is an absolute constant $C_0$
such that for every convex body $K$,
$$
\ell(K^\circ) \le  C_1  \int_0^\infty  \sqrt{ \log N(K, B_2^n,
\e) } \; d \e.
$$

It is possible to slightly improve Dudley's inequality using an
additional volumetric argument. This observation is due to
A.~Pajor.

\begin{lemma}                       \label{dudley}
   There exist absolute constants $C$ and $c$
   such that for a convex body $K$ in $\R^n$
   $$
   \ell(K^\circ)
   \le  C  \int_{c M^*_K}^\infty
              \sqrt{ \log N(K, B_2^n, \e) } \; d \e.
   $$
\end{lemma}

\proof By Dudley's inequality, $\ell(K^\circ) \le  C_1
\int_0^\infty \sqrt{ \log N(K, B_2^n, \e) } \; d \e$. Hence, it
suffices to show that there is some absolute constant $c$ for
which
\begin{equation} \label{eq:integral}
C_1  \int_0^{c M^*_K}
             \sqrt{ \log N(K, B_2^n, \e) } \; d \e
\le  \frac{1}{2} \ell(K^\circ).
\end{equation}
To that end, note that for every $\eps>0$,
\begin{equation}                                \label{N through M}
    N(K, B_2^n, \eps)  \le  \Big(1 +  \frac{2M^*_K}{\e} \Big)^n.
\end{equation}
Indeed, by a standard volumetric argument and Urysohn's
inequality,
\begin{align*}
( N(K, B_2^n, \e) )^{1/n} &\le  \frac{1}{\e}
           \Big( \frac{\vol(K + \e B_2^n)}{\vol(B_2^n)} \Big)^{1/n}
\le  \frac{1}{\e} M^*_{K + \e B_2^n} \\
&\le  \frac{1}{\e} (M^*_K +  M^*_{\e B_2^n}) = \frac{1}{\e} M^*_K
+ 1.
\end{align*}
Thus, by \eqref{N through M}, the integral on the left-hand side
of \eqref{eq:integral} is bounded by
$$C_1 n^{1/2} \int_0^{c M^*_K}
\log^{1/2} (1 + \frac{1}{\e} M^*_K)  \; d \e,
$$
which, after a change of variables, is majorized by
$$
2C_1 n^{1/2}  M^*_K \int_0^{c/2}
  \log^{1/2} (1 + 1 / t)  \; d t \le  C_1 n^{1/2} M^*_K (c /
2)^{1/2} \leq \frac{1}{2} \ell(K^\circ)
$$
for an appropriate choice of $c$.
\endproof

\qquad

\noindent{\bf Proof of Theorem \ref{thm:talagrand}. } By Lemma
\ref{dudley}, there exist absolute constants $C$ and $c$ such that
\begin{equation*}
  E=\ell(K^\circ)
     \le  C\int_{cE / \sqrt{n}}^{\infty}
         \sqrt{\log N (K, B_2^n, t)} \; dt.
\end{equation*}
Since $K \subset  \sqrt{n} B_2^n$, the integrand vanishes for all
$t \ge \sqrt{n}$. Therefore, using Corollary \ref{convex body},
\begin{align*}
E &\le C \int_{cE / \sqrt{n}}^{\sqrt{n}}
            \sqrt{\log N (K, B_2^n, t)} \; dt
= C\sqrt{n}\int_{cE / n}^1 \sqrt{\log N (K, n^{1/2} B_2^n, t)} \; dt \\
&\le C\sqrt{n}\int_{cE / n}^1 \sqrt{VC(K, c t)} \log(2 / t) \; dt,
\end{align*}
as claimed.
\endproof

The main corollary we derive from Theorem \ref{thm:talagrand} is
Elton's Theorem with the optimal dependence on $\d$.

\begin{theorem}                             \label{thm:elton}
   There is an absolute constant $c$ for which the following
   holds.
   Let $x_1, \ldots, x_n$ be vectors in the unit
   ball of a Banach space.
   Assume that for some $\d > 0$
   $$
   \E  \Big\| \sum_{i = 1}^n g_i x_i \Big\|  \ge  \d n.
   $$
   Then there exist two numbers,
   $0 < s < 1$ and $c \d < t < 1$,
   which satisfy that $\sqrt{s} t \log^{2.1} (2 / t)  \ge  \d$,
   and a subset $\s \subset \{ 1, \ldots, n \}$
   of cardinality $|\s|  \ge  s n$,
   such that
   \begin{equation}                      \label{lower l1}
     \Big\| \sum_{i \in \s} a_i x_i \Big\|
     \ge   t \sum_{i \in \s} |a_i|
     \ \ \ \ \text{for all scalars $(a_i)$}.
   \end{equation}
   In particular, we always have $s  \ge  c \d^2$ and $t  \ge  c \d$.
\end{theorem}

\noindent {\bf Proof of Theorem \ref{thm:elton}.} By a
perturbation argument, we may assume that the vectors
$(x_i)_{i=1}^n$ are linearly independent. Hence, using an
appropriate linear transformation we can assume that $X = (\R^n,
\|\cdot\|)$ and that $(x_i)_{i \le n}$ are the unit coordinate
vectors $(e_i)_{i \le n}$ in $\R^n$. Let $K=(B_X)^\circ$ and note
that since $\|e_i\|_X \leq 1$ then $B_1^n \subset K^\circ$.
Therefore, $K \subset B_\infty^n \subset \sqrt{n}B_2^n$.

Let $E=\E\|\sum_{i=1}^n g_i x_i\|_X$. Since $K \subset
B_\infty^n$, then by Theorem \ref{thm:talagrand} there are
absolute constants $c_0$ and $C_0$ such that
\begin{equation*}
\d n \leq E \leq C_0\sqrt{n} \int_{c_0\delta}^{1}
\sqrt{\VC(K,t)}\log(2/t) \; dt.
\end{equation*}
Consider the function
$$
h(t)  =  \frac{c}{t \log^{1.1} (2 / t)}
$$
where the absolute constant $c > 0$ is chosen so that $\int_0^{1}
h(t) \; dt  =  1$. It follows that there exits some $c_0\d  \le  t
\le  1$ such that
$$
\sqrt{\VC(K, c_0t) / n} \cdot \log(2 / t) \ge  \d h(t).
$$
Hence
$$
\VC(K, c_0t)  \ge  \frac{c \d^2}{t^2 \log^{4.2} (2 / t)} n.
$$
Therefore, letting $s = \VC(K, c_0t) / n$ we see that the
announced relation between $s$ and $t$ holds, and that there
exists a subset $\s \subset \{ 1, \ldots, n \}$ of cardinality
$|\s| \ge s n$ such that $(c_0 t / 2) B_\infty^\s  \subset  P_\s
K$. Dualizing, we have $(c_0 t / 2) (K^\circ \cap \R^\s) \subset
B_1^\s$, which completes the proof of the main part of the
theorem.

The ``In particular" part follows trivially.
\endproof

\noindent {\bf Remarks. } Firstly, as the proof shows, the
exponent $2.5$ can be reduced to any number larger than $2$.
Secondly, the relation between $s$ and $t$ in Theorem
\ref{thm:elton} is optimal up to a logarithmic factor for all $0
< \d < 1$. This is seen from by the following example, shown to
us by Mark Rudelson. For $0 < \d < 1 / \sqrt{n}$, the constant
vectors $x_i = \d \sqrt{n} \cdot e_1$ in $X = \R$ show that $s
t^2$ in Theorem \ref{thm:elton} can not exceed $\d^2$. For $1 /
\sqrt{n}  \le \d \le 1$, we consider the body $D = \conv( B_1^n
\cup \frac{1}{\d \sqrt{n}} B_2^n )$ and let $X = (\R^n,
\|\cdot\|_D)$ and $x_i = e_i$, $i = 1, \ldots, n$. Clearly, $\E
\| \sum g_i x_i \|_X  \ge \E \| \sum \e_i e_i \|_D = \d n$. Let
$0 < s, t < 1$ be so that \eqref{lower l1} holds for some subset
$\s \subset \{ 1, \ldots, n \}$ of cardinality $|\s|  \ge  s n$.
This means that $\|x\|_D \ge  t \|x\|_1$ for all $x \in \R^\s$.
Dualizing, we have $\frac{t}{\d \sqrt{n}} \|x\|_2  \le  t
\|x\|_{D^\circ}  \le \|x\|_\infty$ for all $x \in \R^\s$. Testing
this inequality for $x = \sum_{i \in \s} e_i$, we obtain
$\frac{t}{\d \sqrt{n}} \sqrt{|\s|}  \le  1$. This means that $s
t^2  \le  \d^2$.

\qquad

The next application of Theorem \ref{thm:talagrand} is an
improvement of a result of M.~Talagrand \cite{T} which compares
the average over the $\pm$ signs to the minimum over the $\pm$
signs of $\|\sum_{i=1}^n \pm x_i\|$.

\begin{corollary}                       \label{comparison pm}
   Let $x_1, \ldots, x_n$ be vectors in the unit ball
   of a Banach space, and let $M > 0$.
   Fix a number $0  <  \l  <  \log^{-4} (n / M^2)$ and
   assume that
   $$
   \min_{\eta_i = \pm 1} \Big\| \sum_{i \in \s} \eta_i x_i \Big\|
   \le  M |\s|^{1/2}
   \ \ \ \text{for all $\s$ with $|\s|  \le \l n$.}
   $$
   Then
   $$
   \E \Big\| \sum_{i=1}^n g_i x_i \Big\|
   \le C M (n / \l)^{1/2},
   $$
   for some absolute constant $C$.
\end{corollary}

\proof As we did before, we can assume that our Banach space is
$X = (\R^n, \|\cdot\|)$, that $(x_i)_{i=1}^n$ are the unit
coordinate vectors in $\R^n$, and set $K = B_{X^*}$. The
hypothesis of the lemma implies that $\VC(K, M v^{-1/2})  \le v$
if $0  \le  v  \le  \l n$, hence
\begin{equation}                        \label{VC small}
   \VC(K, t)  \le  (M / t)^2
   \ \ \ \text{for $M (\l n)^{-1/2}  \le  t  \le  1$.}
\end{equation}
Let $E = \E \| \sum_{i=1}^n g_i e_i\|_X$. By Theorem
\ref{thm:talagrand}, there are absolute constants $C$ and $c$
such that
$$
E  \le  C \sqrt{n} \int_{c E / n}^1 \sqrt{\VC(K, c t)} \log(2 /
t) \; d t.
$$
If $c E / n  \le  M (\l n)^{-1/2}$, the corollary trivially
follows. Otherwise, if the converse inequality holds, then by
\eqref{VC small},
$$
E  \le  C \sqrt{n} \int_{c E / n}^1 (M / t) \log(2 / t) \; d t
\le  c \sqrt{n} M \cdot \log^2 (n / c E),
$$
and by the assumption on $\lambda$,
$$
E  \le  C \sqrt{n} M \cdot \log^2 (n / M^2) \le  C \sqrt{n} M
\cdot \l^{-1/2},
$$
as claimed.
\endproof

Now we apply Corollary \ref{comparison pm} to compare the type
$2$ constant $T_2(X)$ to the infratype $2$ constant $I_2(X)$ of a
Banach space $X$.

Let $T_2^{(n)} (X)$ and $I_2^{(n)} (X)$ denote the best possible
constants $M$ in \eqref{Tp} and \eqref{Ip}, respectively (with $p
= 2$). So, $T_2^{(n)} (X)$ and $I_2^{(n)} (X)$ measure the
type/infratype $2$ computed on $n$ vectors. Clearly, $I_2 (X)
\le  T_2 (X)$ and $I_2^{(n)} (X)  \le  T_2^{(n)} (X)$.

\begin{corollary}
   Let $X$ be an $n$-dimensional Banach space.
   Then, for every number $0  <  \l  <  \log^{-4} (n / I_2(X)^2)$,
   $$
   T_2 (X)  \le  C \l^{-1/2} \cdot I_2^{(\l n)} (X).
   $$
\end{corollary}
In particular, we obtain
$$
T_2 (X)  \le  I_2 (X)  \cdot  C \log^2 \Big( \frac{n}{I_2(X)^2}
\Big)
   \le  I_2 (X)  \cdot  C \log^2 n.
$$

\proof By \cite{TJ} and \cite{B-K-T} Theorem 3.1, the gaussian
type $2$ can be computed on $n$ vectors of norm one. Precisely,
this means that the constant $T_2 (X)$ equals the smallest
possible constant $M'$ for which the inequality
$$
\E \Big\| \sum_{i=1}^n g_i x_i \Big\| \le M' n^{1/2}
$$
holds for all vectors $x_1, \ldots, x_n$ of norm one. Our
assertion follows from Corollary \ref{comparison pm}.
\endproof

\section{The fat-shattering dimension and covering}
\label{s: empirical measures}

One of the important combinatorial parameters used to measure the
``complexity" of a class of functions is the fat-shattering
dimension, which is a scale-sensitive version of the
Vapnik-Chervonenkis dimension.

\begin{definition} \label{def:fat}
For every $\eps>0$, a set $A=\{ x_1,...,x_n \} \subset \Omega$ is
said to be $\eps$--shattered by $F$ if there is some function
$\g:A \to \R$, such that for every $I \subset \{1,...,n\}$ there
is some $f_I \in F$ for which $f_I(x_i) \geq \g(x_i)+\eps$ if $i
\in I$, and $f_I(x_i) \leq \g(x_i)-\eps$ if $i \not \in I$. Let
\begin{equation*}
{\rm fat}_\eps(F,\Omega)=\sup \Bigl\{|A| \Big| A \subset \Omega, \
A \ {\rm is} \ \eps{\rm -shattered \ by} \ F \Bigr\}.
\end{equation*}
\end{definition}
In cases where the domain is clear, we denote the fat-shattering
dimension of $F$ by ${\rm fat}_\eps(F)$.

If $F$ happens to be a class of Boolean functions, then by
selecting $\g(x_i)=1/2$ we see that ${\rm fat}_{\eps}(F,\Omega)
={\rm VC}(F)$ for every $\eps \leq 1/2$, where ${\rm VC}(F)$ is
the classical Vapnik-Chervonenkis dimension.

Note that the fat-shattering dimension may be controlled by the
generalized VC-dimension, in the following sense. Assume that $F$
is a subset of the unit ball in $L_\infty(\Omega)$, which is
denoted by $B\bigl(L_\infty(\Omega)\bigr)$. Let
$s_n=\{x_1,...,x_n\}$ be a subset of $\Omega$ and set
$F/s_n=\bigl\{\bigl(f(x_1),...,f(x_n)\bigl) \big{|} f \in F
\bigr\} \subset \R^n$. If ${\rm VC}(F/s_n,t) = m $, there is a
subset $\sigma \subset \{1,...,n\}$ of cardinality $m$ such that
$P_\sigma F/s_n \supset \prod_{i \in \sigma} \{a_i,b_i\}$ where
$|b_i-a_i| \geq t$. By selecting $\g(x_i)=(b_i + a_i)/2$ it is
clear that $(x_i)_{i \in \sigma}$ is $t/2$-shattered by $F$, and
thus
$$
\VC(F/s_n,t) \leq {\rm fat}_{t/2}(F,\Omega).
$$

The aim of this section is to bound the entropy of $F$ with
respect to empirical $L_2$ norms. If $s_n=\{x_1,...,x_n\}$ let
$\mu_n$ be the empirical measure supported on $s_n$, that is
$\mu_n=n^{-1}\sum_{i=1}^n \delta_{x_i}$, where $\delta_{x_i}$ is
the point evaluation functional on $x_i$. Empirical covering
numbers play a central role in the theory of empirical processes.
They can be used to characterize classes which satisfy the {\it
uniform law of large numbers} (see \cite{D} or \cite{VW} for a
detailed discussion). It turns out that if $F \subset \bi$ then
$F$ satisfies the uniform law of large numbers with respect to
all probability measures if and only if $\sup_{\mu_n} \log
N\bigl(F,L_2(\mu_n),\e\bigr)=o(n)$ for every $\e>0$, where the
supremum is taken with respect to all empirical measures
supported on at most $n$ elements of $\Omega$. In \cite{A-B-C-H}
it was shown that $F \subset \bi$ satisfies the uniform law of
large numbers if and only if ${\rm fat}_\e(F,\Omega) < \infty$
for every $\e>0$.

Another important application of covering numbers estimates is
the analysis of the {\it uniform central limit property}.
\begin{definition}
Let $F \subset \bi$, set $P$ to be a probability measure on
$\Omega$ and assume $G_P$ to be a gaussian process indexed by $F$,
which has mean $0$ and covariance
\begin{equation*}
\E G_P(f)G_P(g)=\int fgdP-\int fdP \int gdP.
\end{equation*}
A class $F$ is called a universal Donsker class if for any
probability measure $P$ the law $G_P$ is tight in
$\ell_\infty(F)$ and $\nu_n^P=n^{1/2}(P_n -P) \in \ell_\infty(F)$
converges in law to $G_P$ in $\ell_\infty(F)$.
\end{definition}
A property stronger than the universal Donsker property is called
uniform Donsker. For such classes, $\nu_n^P$ converges to $G_P$
uniformly in $P$ in some sense. Instead of presenting the formal
definition of the uniform Donsker property, we mention the
following result of Gin\'{e} and Zinn \cite{GZ}, which
characterizes such classes. Before presenting the result, we
introduce the following notation: for every probability measure
$P$ on $\Omega$, let
$\rho_P^2(f,g)=\E_P(f-g)^2-\bigl(\E_P(f-g)\bigr)^2$, and for
every $\d>0$, set $F_{\d}=\{f-g|f,g \in F, \ \rho_P(f,g) \leq
\d\}$.
\begin{theorem} \cite{GZ}
$F$ is a uniform Donsker property if and only if the following
holds: for every probability measure $P$ on $\Omega$, $G_P$ has a
version with bounded, $\rho_P$-uniformly continuous sample paths,
and for these versions,
\begin{equation*}
\sup_P \E \sup_{f \in F} |G_P(f)| < \infty, \ \ \ \ \lim_{\d \to
0} \sup_P \E \sup_{h \in F_\d} |G_P(h)|=0.
\end{equation*}
\end{theorem}

It is possible to show that the uniform Donsker property is
connected to estimates on covering numbers.
\begin{theorem} \cite{D}
Let $F \subset \bi$. If
\begin{equation*}
  \int_0^\infty  \sup_{n}\sup_{\mu_n}  \sqrt{ \log N
\bigl(F,L_2(\mu_n),\e\bigr) } \;d\eps < \infty,
\end{equation*}
then $F$ is a uniform Donsker class.
\end{theorem}
Having this entropy condition in mind, it is natural to try to
find covering numbers estimates which are ``dimension free", that
is, do not depend on the size of the sample. In the Boolean case,
such bounds where first obtained by Dudley (see \cite{L-T}
Theorem 14.13), and then improved by Haussler \cite{Ha,VW} who
showed that for any empirical measure $\mu_n$ and any Boolean
class $F$,
$$
N(F,L_2(\mu),\e) \leq Cd(4e)^d \eps^{-2d},
$$
where $C$ is an absolute constant and $d={\rm VC}(F)$. In
particular this shows that every VC class is a uniform Donsker
class.

Our goal is to obtain dimension-free estimates on the $L_2$
covering numbers of subsets of $\bi$ using their fat-shattering
dimension, since in many cases it is easier to compute this
parameter than to bound the covering numbers (see, e.g.
\cite{AB}).

Let $F \subset \bi$ and fix a set $s_n \in \Omega$. For every $f
\in F$ let $f/s_n=\sum_{i=1}^n f(x_i)e_i \in F/s_n$. Clearly,
$\|f-g\|_{L_2(\mu_n)}=\|f/s_n -g/s_n\|_{\sqrt{n}B_2^n}$, implying
that for every $t>0$,
\begin{equation} \label{eq:cover}
N\bigl(F,L_2(\mu_n),t\bigr) = N\bigl(F/s_n,\sqrt{n}B_2^n,t\bigr).
\end{equation}

Finally, note that for any $t>0$,
\begin{equation} \label{eq:VCfat}
  {\rm VC}\bigl(F/s_n+\frac{t}{8}B_\infty^n, \frac{t}{2}\bigr) \leq
  {\rm fat}_{\frac{t}{4}}\bigl(F/s_n+\frac{t}{8}B_\infty^n\bigr) \leq
  {\rm fat}_{\frac{t}{8}}(F/s_n) \leq {\rm
  fat}_{\frac{t}{8}}(F).
\end{equation}

\begin{theorem} \label{thm:fat}
There is an absolute constant $C$ such that for any class $F
\subset \bi$, any integer $n$, every empirical measure $\mu_n$ and
every $t>0$,
\begin{equation*}
\log N\bigl(F,L_2(\mu_n),t\bigr) \leq C{\rm fat}_{t/8} (F) \log^2
\frac{2}{t}.
\end{equation*}
\end{theorem}

\proof Let $s_n=\{x_1,...,x_n\}$ be the points on which $\mu_n$
is supported, and apply Corollary \ref{in binfty} for the set
$F/s_n$. We obtain
$$
  \log N (F/s_n, \sqrt{n} B_2^n, t)
   \le  C \log^2 (2 / t )  \cdot  \VC(F/s_n + \frac{t}{8} B_\infty^n, t
/ 2).
$$
Then our claim follows from (\ref{eq:cover}) and (\ref{eq:VCfat}).
\endproof

\begin{remark}
It is possible to show that this bound is essentially tight.
Indeed, fix a class $F \subset \bi$ and put $E(t)=\sup_n
\sup_{\mu_n} \log N\bigl(F,L_2(\mu_n),t \bigr)$ (that is, the
supremum is taken with respect to all the empirical measures
supported on a finite set). By Theorem \ref{thm:fat}, $E(t) \leq
C{\rm fat}_{\frac{t}{8}}(F,\Omega) \log^2 \bigl(2/t \bigr)$. On
the other hand it was shown in \cite{Me} that $E(t) \geq c{\rm
fat}_{16t}(F,\Omega)$ for some absolute constant $c$.
\end{remark}

Comparing the result to Haussler's estimate, one can see that his
bound is recovered up to one logarithmic factor in $1/t$ and the
absolute constant. Indeed, this holds since VC classes satisfy
that ${\rm VC}(F)={\rm fat}_t(F)$ for any $0<t<1/2$.

Now we obtain the following corollary, which extends Dudley's
result from VC classes to the real valued case.

\begin{corollary}
   Let $F  \subset  B(L_\infty(\Omega))$
   and assume that the integral
   $$
   \int_0^1 \sqrt{{\rm fat}_{t/8}(F)} \log \frac{2}{t} \; dt
   $$
   converges.
   Then $F$ is a uniform Donsker class.
\end{corollary}

In particular this shows that if ${\rm fat}_\e (F)$ is ``slightly
better" than $1/\eps^2$, then $F$ is a uniform Donsker class.

\end{document}